\documentclass[twoside, 11pt]{article}

\usepackage{mathrsfs,amsfonts,amsmath}
\usepackage{graphicx}
\usepackage{subfigure}
\usepackage{bm}
\RequirePackage{ifthen,calc}
\usepackage{theorem}
\usepackage[dvips]{color}
\usepackage{booktabs}

\def\E{{\mathbb E}}
\def\P{{\mathbb P}}

 \catcode`@=11
 \def\@evenhead{\hbox to\textwidth{\footnotesize\rm\thepage \hfill
  {\it Single index regression models with randomly left-truncated data}}} % authors name

 \def\@oddhead{\hbox to \textwidth{\footnotesize{\it
 } \hfill\thepage}}% abbreviate title

 \renewcommand{\section}{\makeatletter
 \renewcommand{\@seccntformat}[1]{{\csname the##1\endcsname.}\hspace{0.45em}}
 \makeatother \@startsection
{section}%                                            the name
{1}%                                                  the level
{0pt}%                                                the indent
{\baselineskip}%                                      the beforeskip
{0.5\baselineskip}%                                   the afterskip
{\normalsize\bfseries\mathversion{bold}}}

\renewcommand{\subsection}{\makeatletter
 \renewcommand{\@seccntformat}[1]{{\csname the##1\endcsname.}\hspace{0.45em}}
 \makeatother \@startsection
{subsection}%                                            the name
{1}%                                                  the level
{0pt}%                                                the indent
{\baselineskip}%                                      the beforeskip
{0.5\baselineskip}%                                   the afterskip
{\normalsize\bfseries\mathversion{bold}}}
\catcode`@=12

\newtheorem{theorem}{\noindent Theorem}[section]
\newtheorem{lem}{\noindent Lemma}[section]
\newtheorem{cor}{\noindent Corollary}[section]
\newtheorem{prop}{\noindent Proposition}[section]

\theoremheaderfont{\normalfont\bfseries}
\theorembodyfont{\slshape} {\theorembodyfont{\rmfamily}} {\theorembodyfont{\rmfamily}
}
{\theorembodyfont{\rmfamily} } {\theorembodyfont{\rmfamily}
}
{\theorembodyfont{\rmfamily} } {\theorembodyfont{\rmfamily}
\newtheorem{rem}{\noindent Remark}[section]}
{\theorembodyfont{\rmfamily} } {\theorembodyfont{\rmfamily} } {\theorembodyfont{\rmfamily}
}
 \def\beqlb{\begin{eqnarray}}\def\eeqlb{\end{eqnarray}}
 \def\beqnn{\begin{eqnarray*}}\def\eeqnn{\end{eqnarray*}}
 \numberwithin{equation}{section}
\def\qed{\hfill$\square$\smallskip}

\def\E{{\mathbb{E}}}
\def\P{{\mathbb{P}}}
\begin{document}
\title{\bf  Single Index Regression Models with Randomly Left-truncated Data
\footnotetext{\hspace{-5ex}
${[1]}$ School of Statistics, Shandong University of Finance and Economics,
Jinan 250014,  China. \\
 E-mails: kltgw277519@126.com (Kong, L.), zhangylhappy@163.com (Zhang, Y.), mathdsh@gmail.com (Dai, H.: Corresponding author).
}}
\author{\small Lingtao Kong$^{1}$, Yanli Zhang$^{1}$, Hongshuai Dai$^{1}$}

\maketitle

\begin{abstract}
In this paper, based on the kernel estimator proposed by Ould-Sa\"{i}d and Lemdani (Ann. Instit. Statist. Math. 2006),
we develop some new generalized M-estimator procedures for single index regression models
with left-truncated responses. The consistency and asymptotic normality of our estimators are also established.
Some simulation studies are given to investigate the finite sample performance of the proposed estimators.
\end{abstract}

{\bf Keywords}: Semiparametric regression, single index model, left-truncated data, the product-limit estimator
\vspace{2mm}

\section{Introduction}\label{section1}

In order to avoid the so-called "curse of dimensionality"
in the high dimensional data analysis, many powerful
semiparametric models have been developed to reduce the complexity of
high dimensional data.
One of the popular semiparametric models
is the single index model, which takes the form
\begin{equation}\label{Eq1.1}
Y=g(\theta_0^T\textbf{X};\theta_0)+\epsilon,
\end{equation}
where $Y$ is the response variable, $\textbf{X}\in\mathbb{R}^d (d\geq 2)$
is a covariate vector, $g(\cdot)$ is an unknown univariable measurable
link function, $\epsilon$ is the random error with $\E (\epsilon|\textbf{X})=0$,
$\theta_0\in \mathbb{R}^d$ is the unknown index parameter with
$\|\theta_0\|=1$ (where $\|\cdot\|$ denotes the Euclidean metric)
and the first nonzero component of $\theta_0$ is positive for model identification.
In recent years the single index model has been considered by many authors.
Different methods have been carried out to estimate the index parameter,
 such as the average derivative approach
(Stoker \cite{S1986}, H\"{a}rdle and Tsybakov \cite{HT1993}),
semiparametric least squares estimation (H\"{a}rdle et al. \cite{HHI1993}, Ichimura \cite{I1993}),
semiparametric maximum likelihood estimation (Delecroix et al. \cite{DHP2006}),
 the sliced inverse regression method
 (Duan and Li \cite{DL1991}, Yin and Cook \cite{YC2005}),
 spline estimation (Wang and Yang \cite{WY2009}), and so on.
 Moreover, for the model \eqref{Eq1.1}, Xia et al. \cite{XLTZ2009} considered the goodness-of-fit test.
  Kong and Xia \cite{KX2007} and Wang \cite{W2009} studied the variable selection.
  Xue and Zhu \cite{XZ2006}
  established the empirical likelihood confidence regions for the index parameter, etc.
 Recently, the single index model has been extended to the complex data.
 Bai et al. \cite{BFZ2009} used penalized splines and the method of quadratic
 inference functions to study the single index model for the longitudinal data.
 For the censored data, Lopez \cite{L2009} proposed two semi-parametric M-estimators
 which generalized the estimator of Delecroix, Hristache and Patilea \cite{DHP2006}. Lu and Burke \cite{LB2005} established a $\sqrt{n}$-consistent estimator based on
 the average derivative technique.
 In the case of missing data, Wang et al. \cite{WSHW2010}
 got the estimator of the index parameter and proved the
 asymptotic properties for their estimators.

In practice, the response variable in the model \eqref{Eq1.1} may be left-truncated,
 that is, the variable $(\textbf{X},Y)$ is interfered by another
independent variable $T$ (the truncation variable)
in such a way that we may observe $(\textbf{X},Y)$ and $T$
only when $Y\geq T$, and nothing is observed if $Y<T$.
Truncated data may be encountered in many fields, such as astronomy,
economics, biostatistics and other fields.
Truncated data issues have been investigated extensively
(e.g. Lynden-Bell \cite{L1971}, Woodroofe \cite{W1985}, Stute \cite{S1993},
He and Yang \cite{HY1998a,HY1998b,HY2003},
Stute and Wang \cite{SW2008}, Moreira et al. \cite{MUM2016},
and among others).
Compared with random censored data (or random missing data),
random truncation seems to be more difficult,
since the censored data (or random missing data)
 at least can provide  the information on the censored life-time,
while, in the truncated case, we  observe nothing given $Y<T$.
In this paper, we study the single index model
with left-truncated response.
By extending a kernel estimator for the nonparametric regression
with left-truncated response in Ould-Sa\"{\i}d and Lemdani \cite{OL2006},
 we establish the generalized semiparametric least squares estimators
for the  model \eqref{Eq1.1} in the truncation framework.
The consistency and the asymptotic normality for our estimators are also provided.

The rest of this paper is organized as follows.
In Section 2, we first recall the truncation framework
and  then construct the estimators for $\theta_0$ and
the link function $g(\cdot)$ of the model \eqref{Eq1.1} when
the response is left-truncated.
In Section \ref{section3},
we present the consistency and the asymptotic normality of the estimators.
Section \ref{section4} is devoted to present some simulation studies to test the quality of the estimators
with finite samples. The proofs of our results are collected in Appendix.

\setcounter{equation}{0}
\section{Preliminary}\label{section2}
\subsection{Background for left-truncated data}\label{section2.1}
$(\textbf{X}_j,Y_j,T_j), 1\leq j\leq N$, is a sequence of i.i.d.
random vectors from $(\textbf{X},Y,T)$, where $T$ is the truncation variable.
Throughout this paper, we assume that $T$ is independent of $(\textbf{X},Y)$. Due to
the truncation, we are unable to observe the complete data.
Let
$$(\textbf{X}_{k_i},Y_{k_i},T_{k_i})=:(\textbf{U}_i,V_i,W_i), 1\leq i\leq n,$$
denote the observed sample.
It is obvious that the potential sample size $N$ is unknown
and the observed sample size $n$ is a random variable satisfying $n\leq N$.
We use $\alpha$ to denote the probability we may observe $Y$, that is,
$\alpha=\P(Y\geq T)$. Without loss of generality, we assume $\alpha>0$,
since $\alpha=0$ means that no data can be observed.
For any distribution function $L$,
we use $a_L$ and $b_L$ to stand for the left and right support endpoints of
$L$, respectively.
Define
$$F(y)=\P(Y\leq y),~G(t)=\P(T\leq t),$$ and
$$H(\textbf{x},y)=\P(\textbf{X}\leq \textbf{x},Y\leq y).$$
Let $F^*, G^*$ and $H^{*}$
be the corresponding conditional distributions
of $Y, T$ and $(\textbf{X},Y)$ given $Y\geq T$, respectively,
that is,
\begin{eqnarray*}
F^{*}(v)&=&\alpha^{-1}\int_{-\infty}^{v}G(y)F(dy),\\
G^{*}(v)&=&\alpha^{-1}\int_{-\infty}^{\infty}G(y\wedge v)F(dy),\\
H^{*}(\textbf{u},v)&=&
\alpha^{-1}\int_{-\infty}^{\textbf{u}}\int_{-\infty}^{v}G(y)H(d\textbf{x},dy).
\end{eqnarray*}
It follows from Stute \cite{S1993} and He and Yang \cite{HY2003} that
we can estimate $F^*, G^*$ and $H^{*}$ by $F^{*}_n(v), G^{*}_n(v) $ and $H^{*}_n$, respectively, where
\begin{eqnarray*}
F^{*}_n(v)&=&n^{-1}\sum_{i=1}^nI(V_i\leq v),\\
G^{*}_n(v)&=&n^{-1}\sum_{i=1}^nI(T_i\leq v),\\
H^{*}_n(\textbf{u},v)
&=&n^{-1}\sum_{i=1}^nI(\textbf{U}_i\leq \textbf{u},V_i\leq v).
\end{eqnarray*}

Let $$C(y)=G^{*}(y)-F^*(y)=\alpha^{-1}G(y)[1-F(y)].$$ Then it can be consistently estimated by the following empirical estimator
$$C_n(y)=n^{-1}\sum_{i=1}^nI(W_i\leq y\leq V_i).$$

Next, we introduce  the estimators for  $F(y)$, $G(t)$ and $H(\textbf{x},y)$, respectively. From Lynden-Bell \cite{L1971}, $F$ and $G$ can be estimated  by the so-called Lynden-Bell product-limit  estimators $F_n(y)$ and $G_n(y)$, respectively, where
$$F_n(y)=1-\prod_{V_i\leq y}\Big(1-\frac{1}{nC_n(V_i)}\Big)$$
and
$$G_n(t)=\prod_{W_i> t}\Big(1-\frac{1}{nC_n(W_i)}\Big).$$  On the other hand, He and Yang \cite{HY1998b}
established the following strong consistent estimator $\alpha_n$ for  $\alpha$,
$$\alpha_n=\frac{G_n(y)[1-F_n(y-)]}{C_n(y)},$$
where $F_n(y-)$ denotes the left-continuous version of $F_n(y)$.
Based on the above estimators, He and Yang  \cite{HY2003} got the following nonparametric estimator for
$H(\textbf{x},y)$,
\begin{equation}\label{Eq2.1}
H_n(\textbf{x},y)=\alpha_n\int_{-\infty}^{\textbf{x}}
\int_{-\infty}^{y}\frac{1}{G_n(v)}H_n^{*}(d\textbf{u},dv).
\end{equation}

\subsection{Estimators}

We now come back to our main problem: Estimate the index parameter
$\theta_0$ and the link function $g(\cdot)$ in the model \eqref{Eq1.1} when the
response variable is left-truncated. When data are fully observed, $\theta_0$ and $g(\cdot)$ may be estimated in the following two stages: (i) Estimate the coefficient vector $\theta_0$ ;
(ii) Establish the estimator of the link function $g(\cdot)$
 with the estimator of $\theta_0$ in Step (i). When data are left-truncated, we can still follow the same steps as in the full-data case.
Note that the link function $g$ is unknown.
Due to the left-truncation, we can estimate the link function $g(\theta^T\textbf{u};\theta)$ by the following extended Nadraya-Watson estimator:
\begin{equation}\label{Eq2.2}
\hat{g}_n(\theta^T\textbf{u};\theta)=
\frac{\sum_{i=1}^n V_iG_n^{-1}(V_i)K_{h}(\theta^{T}\textbf{u}-\theta^{T}\textbf{U}_i)}
{\sum_{i=1}^nG_n^{-1}(V_i)K_h(\theta^{T}\textbf{u}-\theta^{T}\textbf{U}_i)},
\end{equation}
where $K_h(\cdot)=K(\cdot/h)$
with $h$ being a bandwidth, and $K(\cdot)$
being a symmetric kernel function with support on $(-1,1)$.
Ould-Sa\"{\i}d and Lemdani \cite{OL2006} and Moreira et al. \cite{MUM2016}
constructed two similar nonparametric estimators
for the regression function with  left-truncated and doubly-truncated responses, respectively.

Similar to the full-data case,  we first estimate $\theta_0$ in the model \eqref{Eq1.1}.
For any measurable function $\varphi(\textbf{u},v)$, under Condition (C1)  (see Section \ref{section3} below),
He and Yang \cite[Theorem 3.2]{HY2003} showed that
$$\int\varphi(\textbf{u},v)H_n(d\textbf{u},dv){\rightarrow}\int\varphi(\textbf{u},v)H(d\textbf{u},dv)~a.s.$$
Hence, we define the estimator $\hat{\theta}_n$ of $\theta_0$ by minimizing $M_n(\theta,\hat{g}_n)$
with
\begin{equation}\label{Eq2.3}
M_n(\theta,\hat{g}_n)=\frac{\alpha_n}{n}\sum_{i=1}^nG_n^{-1}(V_i)
\big[V_i-\hat{g}_n(\theta^{T}\textbf{U}_i;\theta)\big]^2J(\textbf{U}_i),
\end{equation}
where $J(\textbf{u})=I(\textbf{u}\in \mathcal{A})$,
$\mathcal{A}\subset \mathbb{R}^d$, is the trimming function
used to guarantee that the denominator of
$\hat{g}_n(\theta^{T}\textbf{U}_i;\theta)$ is not close to zero.
In the second stage,
 with $\hat{\theta}_n$, the estimator of the link function $g$ is given by
\begin{equation*}
\hat{g}^{*}_n(s;\hat{\theta}_n)=\frac{\sum_{i=1}^n
V_iG_n^{-1}(V_i)K_{h}(s-\hat{\theta}_n^{T}\textbf{U}_i)}
{\sum_{i=1}^nG_n^{-1}(V_i)K_{h}(s-\hat{\theta}_n^{T}\textbf{U}_i)}.
\end{equation*}

\begin{rem}
If there is no truncation, that is, $\alpha=1$, $n=N$ and
$(\textbf{U}_i,V_i)=(\textbf{X}_i,Y_i), 1\leq i\leq N$, then
our estimators  reduce to the ordinary semiparametric least squares estimators.
\end{rem}

\section{Main Results}\label{section3}

In this section, we state the consistency and asymptotic properties of $\hat{\theta}_n$.
We first introduce some notations.
Let $\Theta$ be the set of all unit $d$-vectors with first nonzero component positive.
For any function $f$, let $\nabla_\textbf{x} f$
 (resp. $\nabla^2_{\textbf{x},\textbf{x}}f$)
denote the vector (resp. matrix) of partial derivatives with respect to $\textbf{x}$.

In order to establish our results, we need the following regularity conditions:
\begin{itemize}
\item[(C1)] $F$ and $G$ are continuous with $a_G< a_F$.
\item[(C2)] $\E [Y^2]<\infty.$
\item[(C3)] The set $\mathfrak{X}=Supp(X)$ is a compact subset of $\mathbb{R}^d$.
\item[(C4)] The link function $g(\theta^T\textbf{u})$ is continuous with respect to
 $\theta$ and $\textbf{u}$.
 Furthermore, $g(\theta^T\textbf{u})$ is twice continuously differentiable
with respect to $\theta$, and $\nabla_\theta g$, $\nabla_{\theta,\theta}^2g$ are bounded as functions
of $\theta$ and $\textbf{u}$.\\
\item[(C5)] The kernel function $K$ is a symmetric, positive and twice continuously differentiable
function. Furthermore, $K^{\prime\prime}$ is a Lipschitz continuous function.\\
\item[(C6)] For all $\theta\in\Theta$, the joint density function
$f_{\theta^T\textbf{X},Y}$ of $(\theta^T\textbf{X},Y)$ is twice continuously differentiable with respect to the first variable.\\
\item[(C7)] There exist two Donsker classes $\mathcal{H}_1$ and $\mathcal{H}_2$
such that $$\textbf{u}\rightarrow g(\theta_0^T\textbf{u})\in \mathcal{H}_1
\mbox{ and }  \textbf{u}\rightarrow\nabla_\theta g(\theta_0^T\textbf{u})\in \mathcal{H}_2.$$
\item[(C8)] As $n\rightarrow\infty$, $nh^5(\log n)^{-1}\rightarrow\infty$ and
$n h^7\rightarrow 0$.
\end{itemize}
The continuity of $F$ and $G$ in Condition (C1)
is commonly used in truncated models, see for example
\cite{HY1998a,HY1998b,HY2003,OL2006,Z2011}.
In fact, the continuity  guarantees that
there is no ties in the observed data.
 $a_G<a_F$ in Condition (C1) is needed  for deriving the representation
of $\int \varphi(\textbf{u},v)H_n(d\textbf{u},dv)$
for any measurable function $\varphi$. See
Proposition \ref{propA.1} in Appendix.
Conditions (C2)-(C6) have been widely
used by many authors in the single
index model,  for example, \cite{BL2010,DHP2006,L2009,LPV2013}.
Conditions (C6)-(C8) are used to prove the consistency and the asymptotic normality of our estimator.
The Donsker class in Condition (C7) is also used by \cite{BL2010,L2009,LPV2013}.
For some typical examples of Donsker class,
see van der Vaart and  Wellner \cite[Section 2.10]{VW1996}
and (A.4), (A.5) in \cite{LPV2013}.

\begin{theorem}\label{Th3.1}
Under Conditions (C1)-(C6) and (C8), we have
\begin{equation*}
\hat{\theta}_n\rightarrow\theta_0 \mbox{ in probability}.
\end{equation*}
\end{theorem}

To state the asymptotic normality of $\hat{\theta}_n$,
we first introduce some notations. For any measurable function $\varphi(\textbf{u},v)$, set
\begin{eqnarray}\label{Eq3.1}
\Gamma(\textbf{u},v,\varphi)=\int_{\{v<y\}}[\varphi(\textbf{u},v)-\varphi(\textbf{u},y)]F(dy).
\end{eqnarray}
Moreover, define
\begin{eqnarray}\label{Eq3.2}
\Lambda&=&\E \Big[ \nabla_{\theta}g(\theta_0^T\textbf{U})\nabla_{\theta}
g(\theta_0^T\textbf{U})^TJ(\textbf{U})\Big]
\end{eqnarray}
and
\begin{eqnarray}\label{Eq3.3}
\psi(\textbf{u},v)&=&[v-g(\theta_0^T\textbf{u})]
\nabla_{\theta}g(\theta_0^{T}\textbf{u})J(\textbf{u}).
\end{eqnarray}

\begin{theorem}\label{Th3.2}
Under Conditions (C1) to (C8) , we have
\begin{equation*}
\hat{\theta}_n-\theta_0=n^{-1/2}\Lambda^{-1}W_n+o_\P(n^{-1/2}),
\end{equation*}
where $W_n=n^{-1/2}\sum_{i=1}^n\zeta_i(\psi)$ is a random vector and  for $1\leq i\leq n$,
\begin{equation}\label{Eq3.4}
\zeta_i(\psi)=\frac{\Gamma(\textbf{U}_i,V_i,\psi)}{C(V_i)}
-\int_{T_i}^{V_i}\frac{\Gamma(\textbf{U},v,\psi)}{C^2(v)}F^{*}(dv).
\end{equation}
Hence, as a consequence, we have
$$\hat{\theta}_n-\theta_0\overset{d}{\rightarrow}\mathcal{N}(0,\Lambda^{-1}\Omega \Lambda^{-1})$$
with $$\Omega=Var\Big\{\frac{\Gamma(\textbf{U},V,\psi)}{C(V)}-\int_{T}^{V}
\frac{\Gamma(\textbf{U},v,\psi)}{C^2(v)}F^*(dv)\Big\}.$$
\end{theorem}

The detailed proofs of Theorems \ref{Th3.1} and \ref{Th3.2} are given in Appendix.

\section{Simulation Study}\label{section4}

In this section, we conduct a simulation study to check the finite sample
performance of our estimators.
We conducted the simulation study with the following three different models:

\noindent{{\bf Model 1}:} $Y=-(\theta^T\textbf{X}-1/\sqrt{2})^2+1+\epsilon,$
where $\textbf{X}\sim U[-2,2]\otimes U[-2,2], \epsilon\sim \mathcal{N}(0,0.2^2)$,
the truncated variable $T_1\sim \mathcal{N}(\lambda,1)$ and the true value of the parameter
is
$$\theta_0=(b_1,b_2)^T=\frac{1}{\sqrt{2}}(1,1)^T.$$
This model comes from H\"{a}rdle et al. \cite{HHI1993} and
Lu and Burke \cite{LB2005}.

\noindent{{\bf Model 2}:} $Y=\sin(\theta^T\textbf{X})+\epsilon,$
where $\textbf{X}\sim \mathcal{N}(0,1)\otimes \mathcal{N}(0,1), \epsilon\sim \mathcal{N}(0,0.5^2)$,
the truncated variable $T_2\sim U(-1.5,\lambda)$ and the true value of the parameter
is
$$\theta_0=(b_1,b_2)^T=\frac{1}{\sqrt{5}}(1,2)^T.$$
The second model can be found in Wang et al. \cite{WSHW2010}.

\noindent{{\bf Model 3}:} $Y=\exp\{2\theta^T\textbf{X}\}+\epsilon,$
where $\textbf{X}\sim \mathcal{N}(0,1)\otimes \mathcal{N}(0,1), \epsilon\sim \mathcal{N}(0,1)$,
the truncated variable $T_2\sim \mathcal{N}(\lambda,1)$ and the true value of the parameter
is
$$\theta_0=(b_1,b_2)^T=(0.6,0.8)^T.$$
Moreover, in the above three models, the variables
$\textbf{X}, \epsilon$ and $T$ are mutually independent.

Here, we should point out that, from Section \ref{section2.1}, the estimators $F_n$ and $G_n$
depend highly on the behavior of the estimator $C_n$,
while $C_n$ may be zero with the truncated data.
To overcome this problem, similar to Woodroofe \cite{W1985} and Zhou \cite{Z2011}, in the simulation study,
we replaced $C_n(y)$ by
$$\tilde{C}_n(y)=\max\Big\{C_n(y),\frac{1}{n}+\frac{1}{n^2}\Big\},~\mbox{for any}~y \in(V_{(1)},V_{(n)}),$$
where $V_{(1)}, V_{(n)}$ are the ordered statistics.
Moreover,  Stute and Wang \cite{SW2008} proved that the corresponding
estimators based on $C_n(y)$ and $\tilde{C}_n(y)$ are asymptotically equivalent at the $\sqrt{n}$-rate.

In our simulations, we assumed  that the complete data size $N$
is fixed and the observed data size $n$ is random for convenience (you may also
set $n$ be fixed and $N$ be random).
For each model, we performed $500$ repetitions for  each setting $(N,\alpha)$, where
the sample size $N\in\{50, 100, 200\}$
and the proportions of truncated data $1-\alpha=\P(Y<T)\in\{10\%, 20\%, 40\%\}$.
We chose the Epanechnikov kernel function
$K(u)=\frac{3}{4}(1-u^2)I\{|u|\leq 1\}$ and the bandwidth sequence $h=n^{-1/5}(\log n)^{1/5}$  to
compute $\hat{\theta}_n$.  The bias and the mean squared error (MSE)
for $\hat{\theta}_n$ were computed. The corresponding results are presented  in Tables \ref{table1}-\ref{table3}.

\begin{table}[htbp]
\caption{Simulation results for Model 1.}
\begin{tabular*}{15cm}{@{\extracolsep{\fill}}lllllll}
\toprule
$\lambda$&$1-\alpha$&N&\multicolumn{2}{c}{Bias}&\multicolumn{2}{c}{MSE}\\
\cline{4-5}\cline{6-7}
&&&$b_1$&$b_2$&$b_1$&$b_2$\\
\midrule
  -0.72&0.4&50&$~~7.38\times10^{-3}$&$-4.01\times10^{-3}$&$3.20\times10^{-3}$&$1.57\times10^{-3}$\\
  &&100&$~~1.02\times10^{-3}$&$~~3.58\times10^{-5}$&$6.92\times10^{-4}$&$7.10\times10^{-4}$\\
  &&200&$-3.95\times10^{-4}$&$~~7.83\times10^{-4}$&$2.75\times10^{-4}$&$2.73\times10^{-4}$\\
  -2.4&0.2&50&$~~1.63\times10^{-3}$&$-8.39\times10^{-4}$&$5.82\times10^{-4}$&$5.48\times10^{-4}$\\
  &&100&$~~3.42\times10^{-4}$&$-2.63\times10^{-5}$&$2.27\times10^{-4}$&$2.19\times10^{-4}$\\
  &&200&$-7.93\times10^{-4}$&$~~9.21\times10^{-4}$&$9.01\times10^{-5}$&$9.20\times10^{-5}$\\
  -3.5&0.1&50&$-2.20\times10^{-4}$&$~~8.69\times10^{-4}$&$4.53\times10^{-4}$&$4.65\times10^{-4}$\\
  &&100&$~~1.42\times10^{-3}$&$-1.26\times10^{-3}$&$1.60\times10^{-4}$&$1.57\times10^{-4}$\\
  &&200&$~~3.10\times10^{-4}$&$~~4.06\times10^{-4}$&$6.74\times10^{-5}$&$6.84\times10^{-5}$\\
\bottomrule
\end{tabular*}
\label{table1}
\end{table}

\begin{table}[htbp]
\caption{Simulation results for Model 2.}
\begin{tabular*}{15cm}{@{\extracolsep{\fill}}lllllll}
\toprule
$\lambda$&$1-\alpha$&N&\multicolumn{2}{c}{Bias}&\multicolumn{2}{c}{MSE}\\
\cline{4-5}\cline{6-7}
&&&$b_1$&$b_2$&$b_1$&$b_2$\\
\midrule
  0.92&0.4&50&$-2.29\times10^{-2}$&$-8.76\times10^{-3}$&$2.87\times10^{-2}$&$7.41\times10^{-3}$\\
  &&100&$-1.50\times10^{-2}$&$-3.30\times10^{-3}$&$1.52\times10^{-2}$&$4.11\times10^{-3}$\\
  &&200&$-4.35\times10^{-3}$&$-2.25\times10^{-3}$&$6.36\times10^{-3}$&$1.57\times10^{-3}$\\
  -0.13&0.2&50&$-1.47\times10^{-2}$&$-6.33\times10^{-2}$&$1.93\times10^{-2}$&$5.16\times10^{-3}$\\
  &&100&$-3.29\times10^{-3}$&$-5.70\times10^{-3}$&$1.03\times10^{-2}$&$2.85\times10^{-3}$\\
  &&200&$-5.98\times10^{-3}$&$~~7.52\times10^{-5}$&$4.20\times10^{-3}$&$1.01\times10^{-3}$\\
  -0.75&0.1&50&$-1.19\times10^{-2}$&$-5.10\times10^{-3}$&$1.60\times10^{-2}$&$3.82\times10^{-3}$\\
  &&100&$-3.97\times10^{-3}$&$-3.41\times10^{-3}$&$7.60\times10^{-3}$&$2.07\times10^{-3}$\\
  &&200&$-5.60\times10^{-3}$&$~~2.05\times10^{-4}$&$3.74\times10^{-3}$&$9.06\times10^{-4}$\\
\bottomrule
\end{tabular*}
\label{table2}
\end{table}

From Tables \ref{table1}-\ref{table3}, we can see that our estimator
$\hat{\theta}_n$ performs well.
Moreover, the performance of  $\hat{\theta}_n$  become better and better as the sample size $N$ increases.
We also observe  that the quality of our estimator in each model is slightly
affected by the proportion of the truncated data, $1-\alpha$,
and shrinks as the proportion becomes larger.

\begin{table}[htbp]
\caption{Simulation results for Model 3.}
\begin{tabular*}{15cm}{@{\extracolsep{\fill}}lllllll}
\toprule
$\lambda$&$1-\alpha$&N&\multicolumn{2}{c}{Bias}&\multicolumn{2}{c}{MSE}\\
\cline{4-5}\cline{6-7}
&&&$b_1$&$b_2$&$b_1$&$b_2$\\
\midrule
  0.97&0.4&50&$-7.38\times10^{-3}$&$~~7.81\times10^{-5}$&$5.52\times10^{-3}$&$3.20\times10^{-3}$\\
  &&100&$-7.08\times10^{-3}$&$~~2.47\times10^{-3}$&$3.11\times10^{-3}$&$1.42\times10^{-3}$\\
  &&200&$-1.35\times10^{-3}$&$-1.22\times10^{-3}$&$2.34\times10^{-3}$&$1.23\times10^{-3}$\\
  -0.20&0.2&50&$-6.83\times10^{-3}$&$~~6.10\times10^{-5}$&$5.11\times10^{-3}$&$2.99\times10^{-3}$\\
  &&100&$-2.39\times10^{-3}$&$-3.18\times10^{-3}$&$3.08\times10^{-3}$&$2.30\times10^{-3}$\\
  &&200&$-7.97\times10^{-4}$&$-1.47\times10^{-3}$&$2.25\times10^{-3}$&$1.06\times10^{-3}$\\
  -4.3&0.1&50&$~~3.79\times10^{-5}$&$~~5.33\times10^{-4}$&$4.01\times10^{-4}$&$4.05\times10^{-4}$\\
  &&100&$~~5.38\times10^{-4}$&$-3.61\times10^{-4}$&$1.25\times10^{-4}$&$1.24\times10^{-4}$\\
  &&200&$-1.87\times10^{-4}$&$~~2.77\times10^{-4}$&$6.27\times10^{-5}$&$6.41\times10^{-5}$\\
\bottomrule
\end{tabular*}
\label{table3}
\end{table}

\begin{figure}[htb]
    \centering
\subfigure{\includegraphics[scale=0.3]{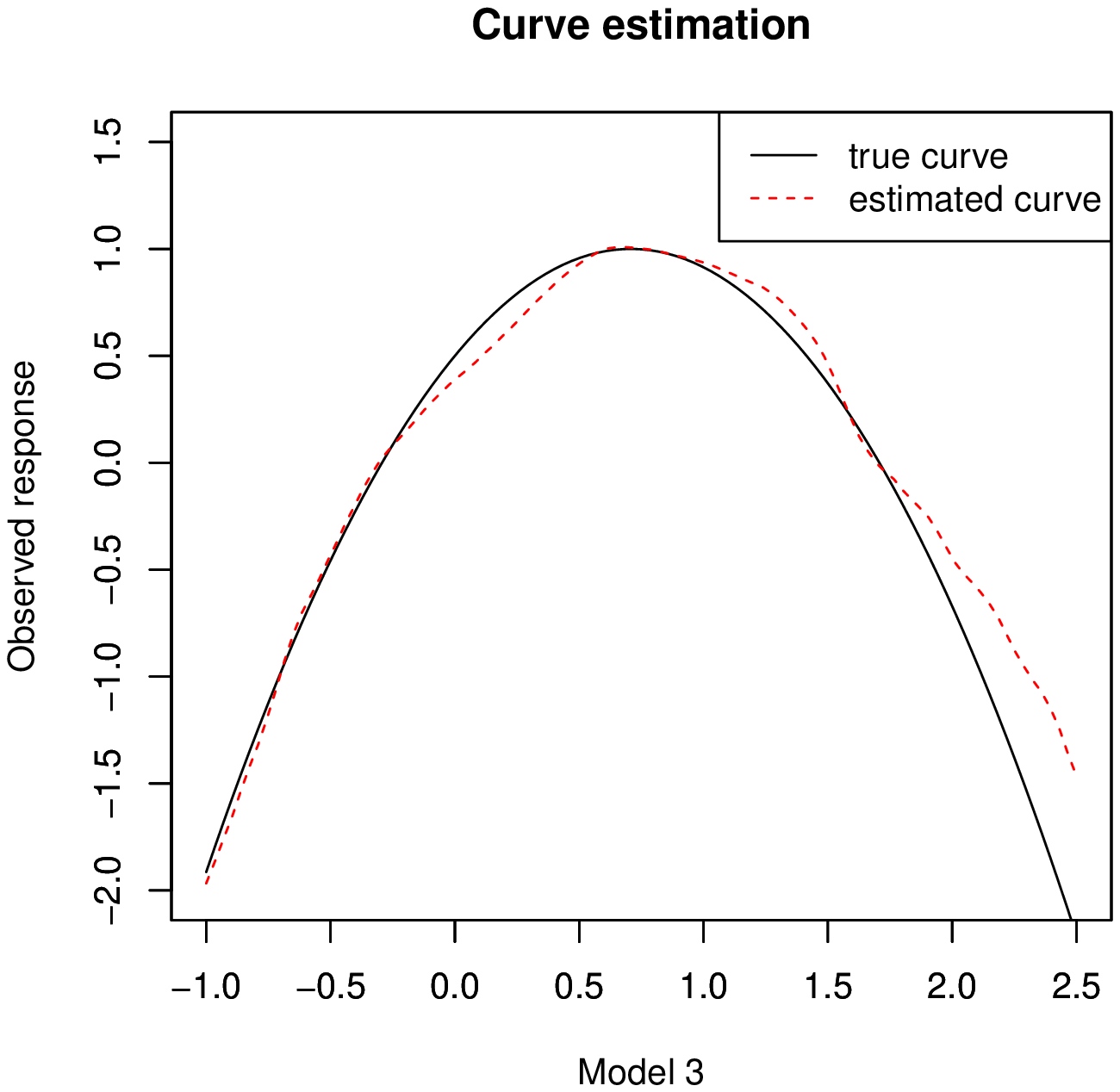}}
\subfigure{\includegraphics[scale=0.3]{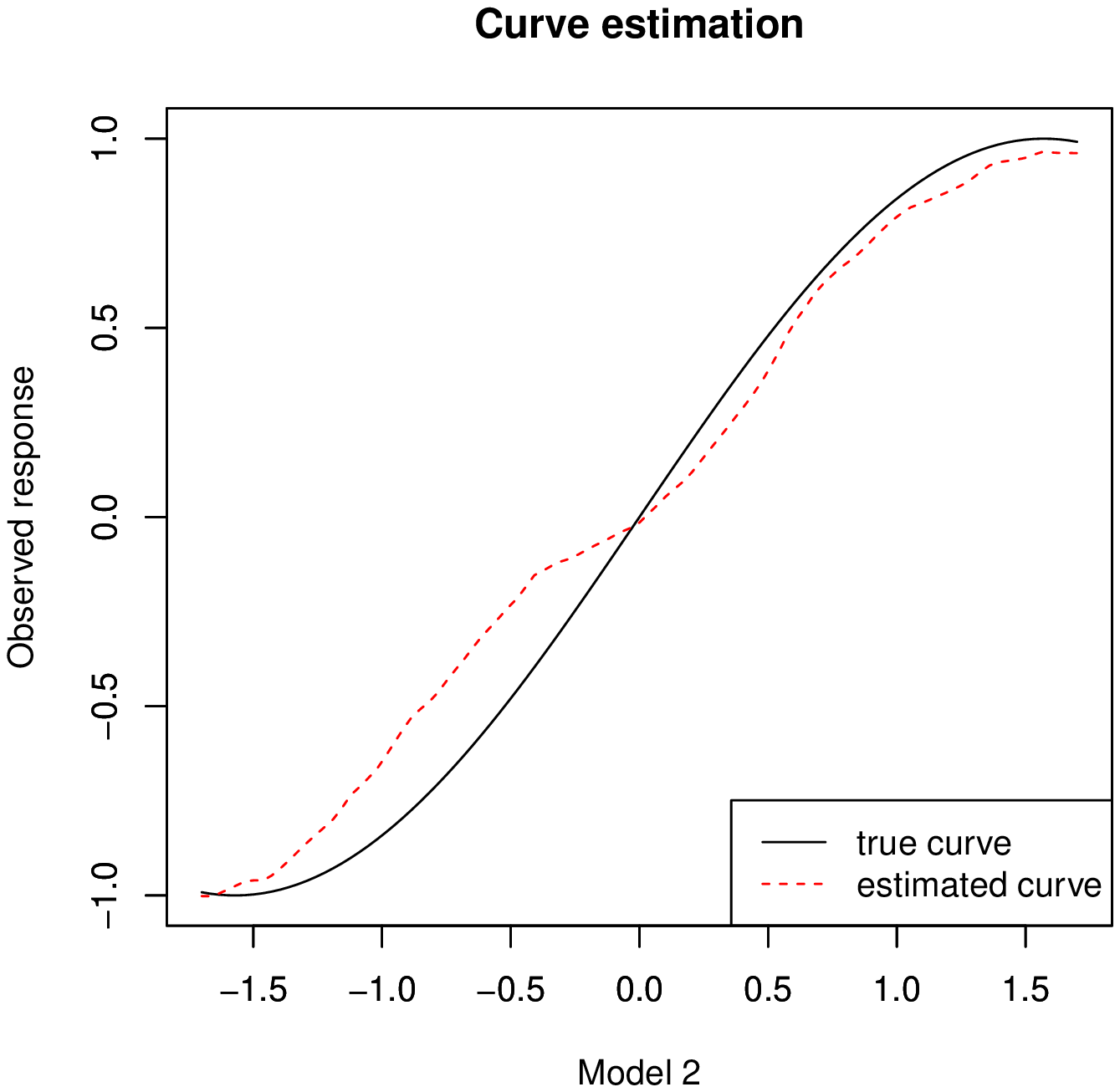}}
\subfigure{\includegraphics[scale=0.3]{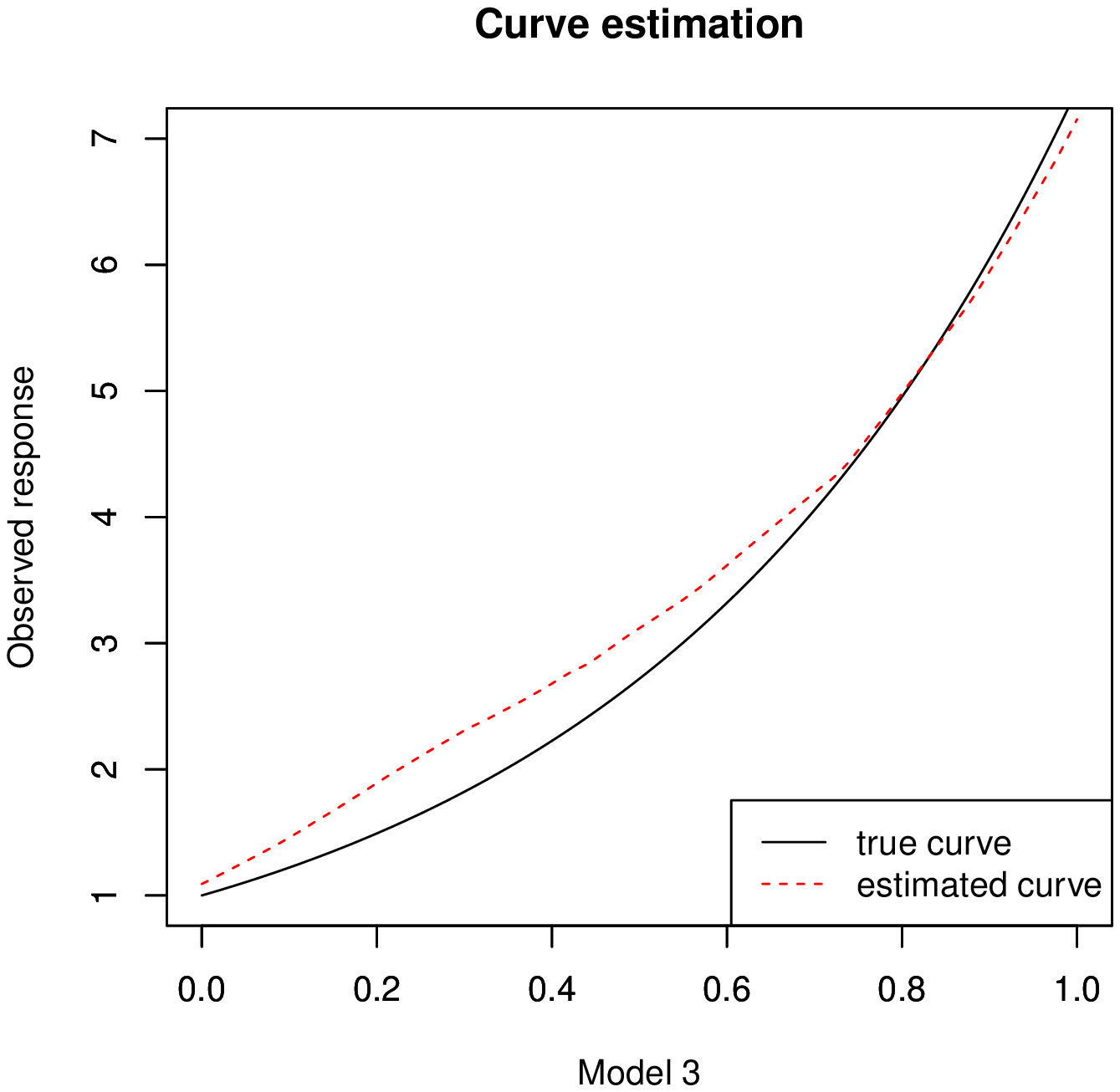}}
\caption{Curve estimations for Models 1$\sim$3 with $N=200$.}
\label{figure1}
\end{figure}

\begin{figure}[htb]
    \centering
\subfigure{\includegraphics[scale=0.3]{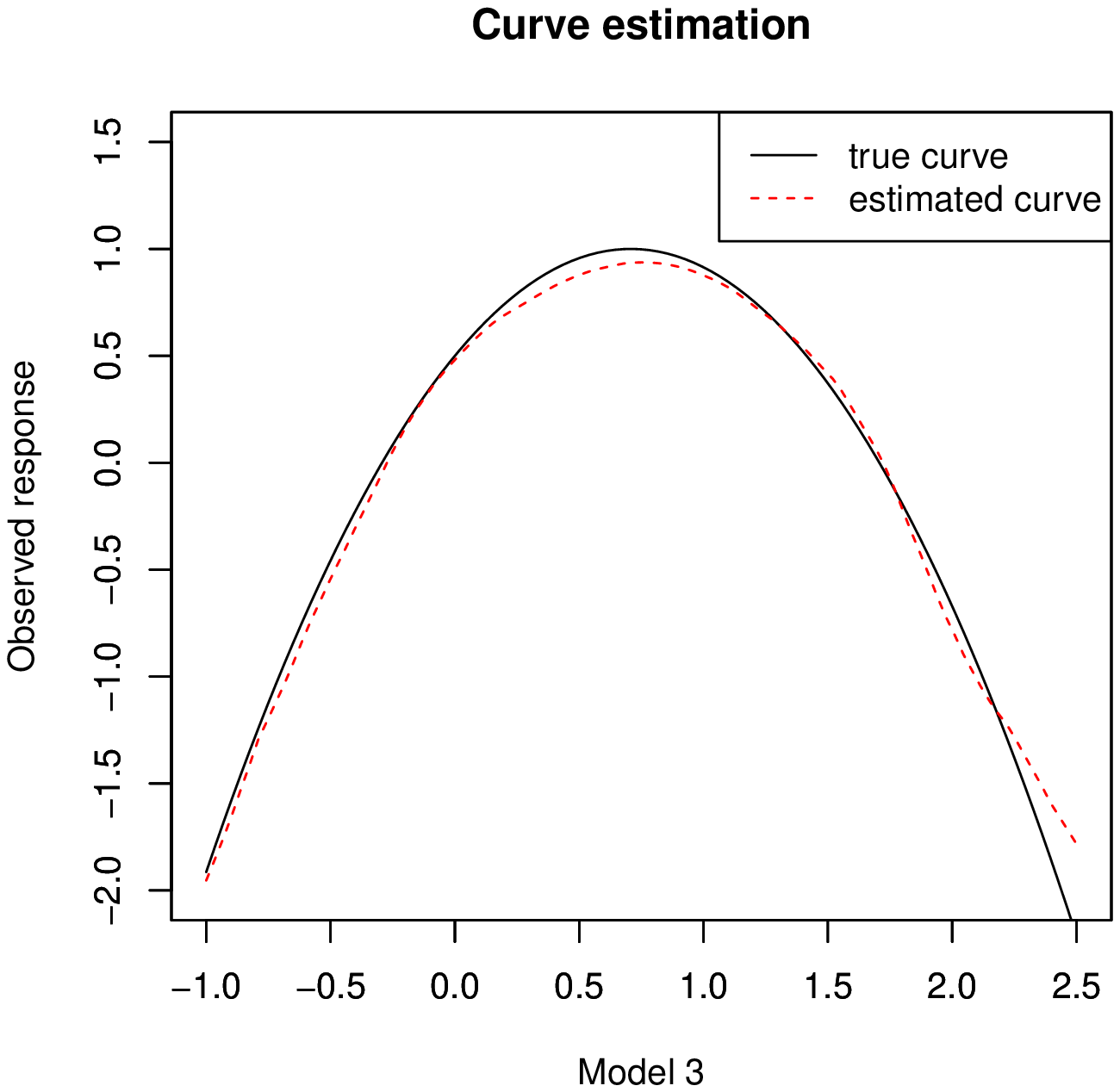}}
\subfigure{\includegraphics[scale=0.25]{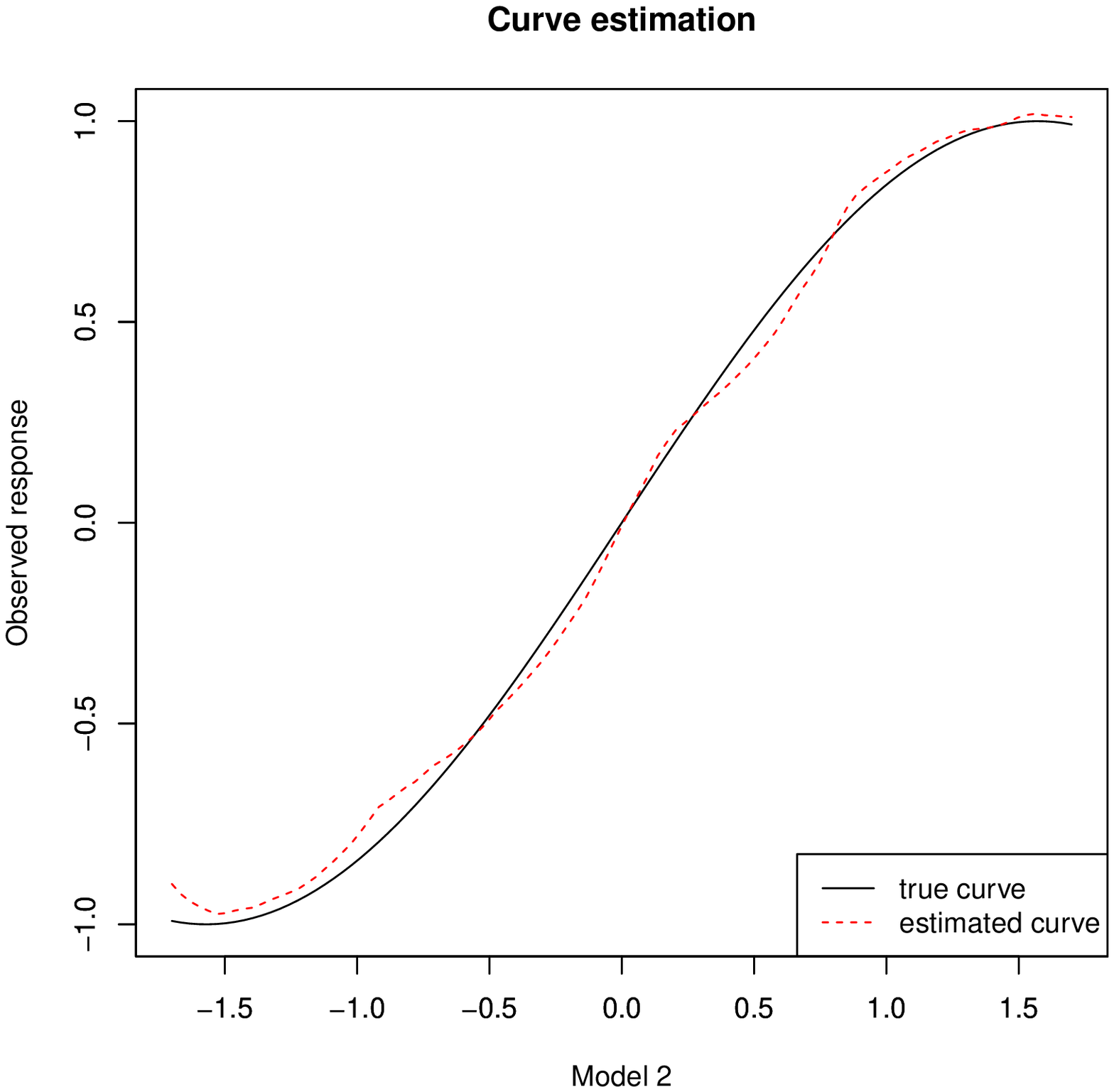}}
\subfigure{\includegraphics[scale=0.3]{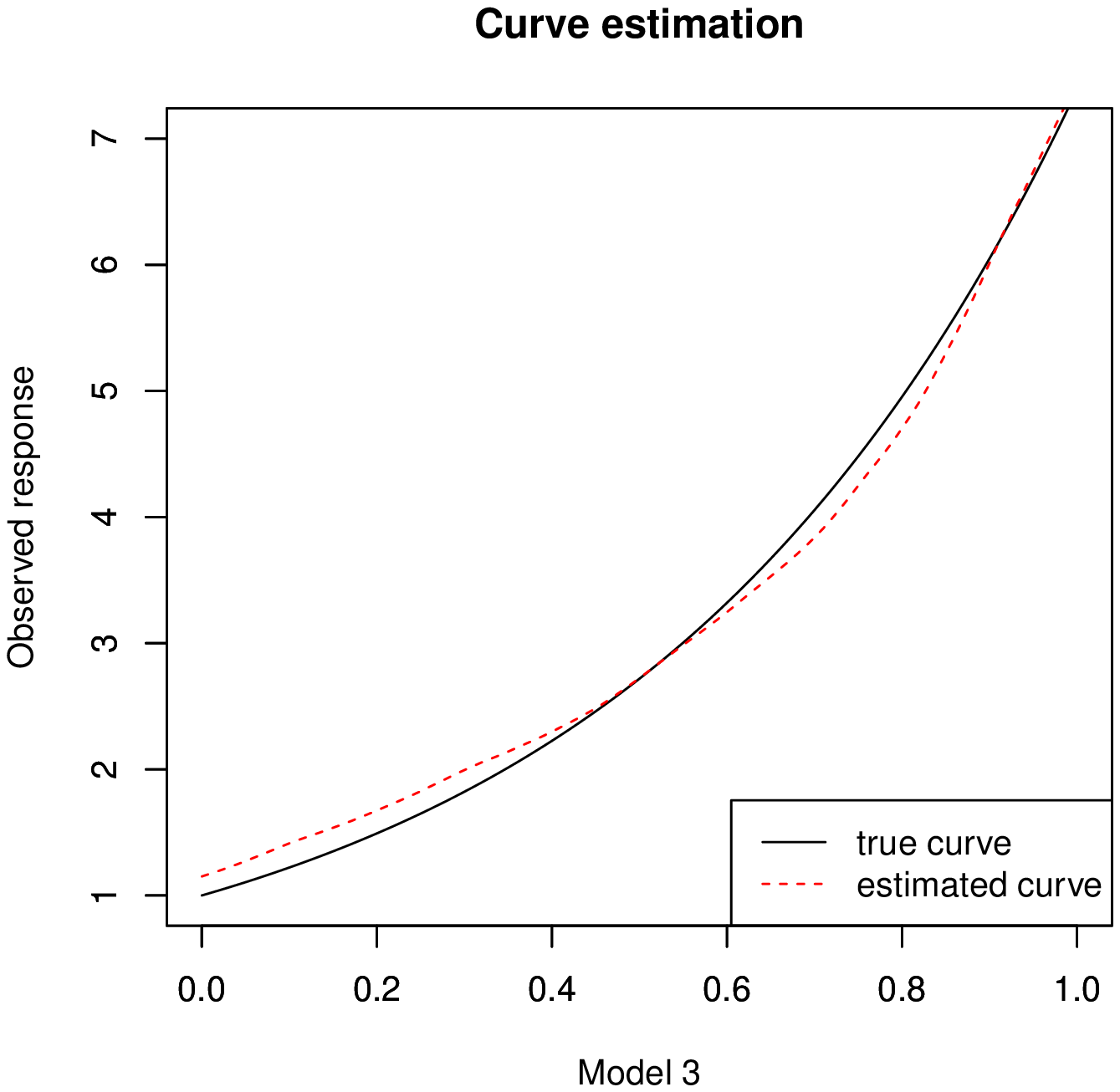}}
\caption{Curve estimations for Models 1$\sim$3 with $N=500$.}
\label{figure2}
\end{figure}

Corresponding to $N=200$ and $N=500$, the curves of $\hat{g}_n$ for
three models  with $1-\alpha=20\%$ are graphed in Figures 1 and 2, respectively.   The appearance of the estimated curves is very similar to that of the true curves. Figures \ref{figure1} and \ref{figure2} suggest that our estimators work well too.

\section{Concluding Remarks}

In this paper, we have considered the single index models
under random truncated framework.
The estimators of the index parameter
and the link function are established based on the kernel
estimator proposed by \cite{OL2006}.
 Our estimators possess the consistency and the asymptotic normality.
Simulations indicate that  the proposed method performs well.

Most of statistical methods dealing with the truncated data, including
the Lynden-Bell estimator \cite{L1971}, rely heavily on the quasi-independence (independence)
between the truncated random $T$ and the interest variable $Y$. See, for example,
\cite{HY1998a,HY1998b,HY2003,OL2006,SW2008,Z2011}.
However, the quasi-independence (independence) may fail in many situations.
For example, Chaieb et al. \cite{CRA2006} introduced a copular dependency between $T$ and $Y$,
and established some modified estimators for the distribution
functions.
Thus, it will be interesting  to extend some similar ideas to our setting.
We will investigate this in the future.

\begin{appendix}
\section{Appendix: Technical Proofs}

\subsection{Representation of $\int\varphi(\textbf{u},v)H_n(d\textbf{u},dv)$}
In this subsection, we
study the representation of
$\int \varphi(\textbf{u},v)H_n(d\textbf{u},dv)$
for any measurable function $\varphi(\textbf{u},v)$, which  is essential for the proof of
Theorem \ref{Th3.2}.
 \begin{prop}\label{propA.1}
Let $\varphi(\textbf{u},v)$ be any measurable function satisfying
\begin{equation}\label{EqA.1}
\int\frac{\varphi^2(\textbf{u},v)}{G(v)} H(d\textbf{u},dv)<\infty.
\end{equation}
Then, under Condition (C1), we have
\begin{eqnarray*}
&&\int\varphi(\textbf{u},v)[H_n(d\textbf{u},dv)-H(d\textbf{u},dv)]\\
&&\quad=
\int\frac{\Gamma(\textbf{u},v,\varphi)}{C(v)}[H_n^{*}(d\textbf{u},dv)-H^{*}(d\textbf{u},dv)]\\
&&\quad\quad-\int\frac{C_n(v)-C(v)}{C^2(v)}\Gamma(\textbf{u},v,\varphi)H^{*}(d\textbf{u},dv)+o_{\P}(n^{-\frac{1}{2}}),
\end{eqnarray*}
where $\Gamma(\textbf{u},v,\varphi)$ is defined by \eqref{Eq3.1}.
\end{prop}
\begin{rem}
Proposition \ref{propA.1} extends Theorem 1.1 in Stute and Wang \cite{SW2008}
which dealt with the representation when the covariables are absent.
It is obvious that the inclusion of covariables will enlarge
the class of possible applications.
\end{rem}
\noindent{\it Proof:}
The proof of Proposition \ref{propA.1} is similar to that of
Stute and Wang \cite[Theorem 1.1]{SW2008}
which studies the representation without covariables. The rest of the proof is devoted to some modifications.

We first introduce an asymptotically equivalent estimator $\hat{H}_n$ of $H_n$,
which is defined by
\begin{equation*}
\int\varphi(\textbf{u},v)\hat{H}_n(d\textbf{u},dv):=
\int \frac{\varphi(\textbf{u},v)\lambda_n(v)}{C_n(v)}H_n^*(d\textbf{u},dv),
\end{equation*}
where
$$\lambda_n(v)=\exp\Big\{n\int_{-\infty}^{v-}\log\big[1-\frac{1}{1+nC_n(y)}\big]F_n^{*}(dy)\Big\}$$
with $\int_{-\infty}^{v-}$ denoting the integral on the interval $(-\infty,v)$.
Set $$\lambda(v)=1-F(v)=\exp\Big\{-\int_{-\infty}^v\frac{F^*(dy)}{C(y)}\Big\}.$$
Similar to the proofs of Lemma 3.2 and Corollaries 3.1-3.3 in Stute and Wang \cite{SW2008},
we  obtain that
\begin{eqnarray}\label{EqA.2}
\int\varphi(\textbf{u},v)\hat{H}_n(d\textbf{u},dv)=L_{n1}+L_{n2}+L_{n3}+o_{\P}(n^{-1/2})
\end{eqnarray}
with
\begin{eqnarray*}
&&L_{n1}=\int\frac{\varphi(\textbf{u},v)\lambda(v)}{C(v)}H_n^*(d\textbf{u},dv)
+\int\frac{\varphi(\textbf{u},v)\lambda(v)(C(v)-C_n(v))}{C^2(v)}H_n^*(d\textbf{u},dv),\\
&&L_{n2}=-\int\frac{\varphi(\textbf{u},v)\lambda(v)}{C(v)}
\int_{-\infty}^{v-}\frac{F^*_n(dy)-F^*(dy)}{C(y)}H_n^*(d\textbf{u},dv),\\
&&L_{n3}=\int\frac{\varphi(\textbf{u},v)\lambda(v)}{C(v)}
\int_{-\infty}^{v-}\frac{C_n(y)-C(y)}{C^2(y)}F_n^*(dy)H_n^*(d\textbf{u},dv).
\end{eqnarray*}
It follows from \eqref{EqA.2} and Theorems 5.3.2 and 5.3.3 in Serfling \cite{S1980} that
Proposition \ref{propA.1} holds for $\hat{H}_n$ instead of $H_n$.
By applying the SLLN for U-statistics,
we  obtain that Proposition \ref{propA.1} also holds for $H_n$. %For more details, see \cite[P. 621-622]{SW2008}.
\qed

From Proposition \ref{propA.1}, we get an i.i.d. representation of
$\int\varphi(\textbf{u},v) H_n(d\textbf{u},dv)$
and the asymptotic normality. For similar results of the censored data, refer to Theorem 1.1 in Stute \cite{S1996}.

\begin{cor}\label{corA.1}
Under the assumptions of Proposition \ref{propA.1}, we have
$$\sqrt{n}\int\varphi(\textbf{u},v)[H_n(d\textbf{u},dv)-H(d\textbf{u},dv)]=
n^{-\frac{1}{2}}\sum_{i=1}^n\zeta_i(\varphi)+o_{\P}(1),$$
where $\zeta_i(\cdot)$ is defined by \eqref{Eq3.4}.
\end{cor}

\begin{cor}\label{corA.2}
Under the assumptions of Proposition \ref{propA.1}, we have
$$\sqrt{n}\int\varphi(\textbf{u},v)[H_n(d\textbf{u},dv)-H(d\textbf{u},dv)]
\overset{d}{\rightarrow}\mathcal{N}(0,\sigma^2),$$
where
\begin{equation*}
\sigma^2=Var\Big\{\frac{\Gamma(\textbf{U},V,\varphi)}{C(V)}-\int_{T}^{V}
\frac{\Gamma(\textbf{U},v,\varphi)}{C^2(v)}F^*(dv)\Big\},
\end{equation*}
and the function $\Gamma$ is defined by \eqref{Eq3.1}.
\end{cor}

\begin{rem}\label{remA.1}
$a_G<a_F$ in Condition (C1), together with
\begin{equation}\label{EqA.3}
\int \varphi^2(\textbf{u},v)H(d\textbf{u},dv)<\infty,
\end{equation}
means that \eqref{EqA.1} holds.
Note that \eqref{EqA.3}
is the standard moment condition when the data are complete.
\end{rem}

\begin{rem}\label{remA.2}
S\'{a}nchez Sellero et al. \cite[Theorem 1]{S2006} introduced an i.i.d. representation for
 the product-limit integrals under truncation and censoring with covariables.
However, Theorem 1 in  \cite{S2006} requires that  the following two integrals
\begin{eqnarray}\label{EqA.4}
\int \Phi(\textbf{u},v)[1-F(v)]^{-5}H(d\textbf{u},dv) \mbox{ and }
\int \Phi^2(\textbf{u},v)[1-F(v)]^{-3}H(d\textbf{u},dv)
\end{eqnarray}
are finite where $\Phi$ is an envelope for the class
$\{\varphi(\textbf{u},v)\}$
 (we refer to var der Vaart and Wellner \cite[P. 84]{VW1996} for the definition of the envelope function).
It is obvious that both the integrals in \eqref{EqA.4}
equal infinity when $\Phi$ is a finite constant. Hence, the representation
proposed by \cite{S2006} can only be applied to the functions
converging to zero as $\textbf{u},v\rightarrow\infty$.
However, Proposition \ref{propA.1} only requires a finite second moment condition.
\end{rem}

\subsection{Difference between $g$ and $\hat{g}_n$}

To prove the consistency and the asymptotic normality of $\hat{\theta}_n$,
we need to study the difference between  $g$ and $\hat{g}_n$ defined by \eqref{Eq2.2}.

Let $$\phi_{\theta^T\textbf{X}}(s)=\int yf_{\theta^T\textbf{X},Y}(s,y)dy,$$
where $f_{\theta^T\textbf{X},Y}(s,y)$
is the joint density function of $(\theta^T\textbf{X},Y)$.
Hence, we can rewrite the link function $g(\cdot)$ in \eqref{Eq1.1} as
$$g(s;\theta_0)=\frac{\phi_{\theta_0^T\textbf{X}}(s)}{f_{\theta_0^T\textbf{X}}(s)},$$
where $f_{\theta^T\textbf{X}}(\cdot)$ is the density function of $\theta^T\textbf{X}$.
Define $$\hat{f}_{\theta^T\textbf{X},n}(s)=
\frac{\alpha_n}{nh}\sum_{i=1}^n\frac{1}{G_n(V_i)}
K\Big(\frac{s-\theta^T\textbf{U}_i}{h}\Big),$$
$$\hat{\phi}_{\theta^T\textbf{X},n}(s)=
\frac{\alpha_n}{nh}\sum_{i=1}^n\frac{V_i}{G_n(V_i)}
K\Big(\frac{s-\theta^T\textbf{U}_i}{h}\Big).$$
Note that
$$\hat{g}_n(s;\theta)=\frac{\hat{\phi}_{\theta^T\textbf{X},n}(s)}{\hat{f}_{\theta^T\textbf{X},n}(s)}.$$
The following lemmas \ref{lemA.1} to \ref{lemA.5} study the distance between $g$
and $\hat{g}_n$. Before we state them, we first introduce two equivalent estimators. Set
$$\tilde{f}_{\theta^T\textbf{X},n}(s)=
\frac{\alpha}{nh}\sum_{i=1}^n\frac{1}{G(V_i)}
K\Big(\frac{s-\theta^T\textbf{U}_i}{h}\Big),$$
$$\tilde{\phi}_{\theta^T\textbf{X},n}(s)=\frac{\alpha}{nh}\sum_{i=1}^n\frac{V_i}{G(V_i)}
K\Big(\frac{s-\theta^T\textbf{U}_i}{h}\Big).$$

\begin{lem}\label{lemA.1}
Under Conditions (C1), (C2), (C5) and (C8), we have
\begin{eqnarray}\label{EqA.5}
&&\sup_{\textbf{u}\in\mathfrak{X},\theta\in\Theta}\Big|\tilde{f}_{\theta^T\textbf{X},n}(\theta^T\textbf{u})
-\E\tilde{f}_{\theta^TX,n}(\theta^T\textbf{u})\Big|=
O\Big(\sqrt{\frac{\log n}{n h}}\Big)~ a.s.
\end{eqnarray}
and
\begin{eqnarray}\label{EqA.6}
&&\sup_{\textbf{u}\in\mathfrak{X},\theta\in\Theta}\Big|\tilde{\phi}_{\theta^T\textbf{X},n}(\theta^T\textbf{u})
-\E\tilde{\phi}_{\theta^TX,n}(\theta^T\textbf{u})\Big|=
O\Big(\sqrt{\frac{\log n}{n h}}\Big)~a.s.
\end{eqnarray}
\end{lem}
\noindent{\it Proof:} Let
$\mathfrak{F}_1=\{\frac{1}{G(y)}\}$.
From Gin\'{e} and Guillou \cite[Lemma 3(a)]{GG1999},
under Condition (C1)
the class $\mathfrak{F}_1$
is a V-C subgraph class (see Gin\'{e} and Guillou \cite[P. 2049]{GG1999})
with the envelope $\frac{1}{G(a_F)}.$
Hence, the assumptions of Theorem 1
 in Einmahl and Mason \cite{EM2000} hold under Conditions (C1), C(5) and (C8).
 Thus, by applying Theorem 1 in
 \cite{EM2000}, we  conclude \eqref{EqA.5}.

The proof of \eqref{EqA.6} is similar to that of \eqref{EqA.5},
but using the V-C subgraph class
 $\mathfrak{F}_2=\{\frac{y}{G(y)}\}$ with the envelope $\frac{y}{G(a_F)}$
 instead of the class $\mathfrak{F}_1$.
\qed

\begin{lem}\label{lemA.2}
Under Conditions (C1), (C2), (C5), (C6) and (C8), we have
\begin{eqnarray*}
&&\sup_{\textbf{u}\in\mathfrak{X},\theta\in\Theta}
\Big|f_{\theta^T\textbf{X}}(\theta^T\textbf{u})
-\tilde{f}_{\theta^T\textbf{X},n}(\theta^T\textbf{u})\Big|=
O\Big(\sqrt{\frac{\log n}{nh}}+h^2\Big)~a.s.
\end{eqnarray*}
and
\begin{eqnarray*}
&&\sup_{\textbf{u}\in\mathfrak{X},\theta\in\Theta}
\Big|\phi_{\theta^T\textbf{X}}(\theta^T\textbf{u})
-\tilde{\phi}_{\theta^T\textbf{X},n}(\theta^T\textbf{u})\Big|=
O\Big(\sqrt{\frac{\log n}{nh}}+h^2\Big)~a.s.
\end{eqnarray*}
\end{lem}
\noindent{\it Proof:}
From Lemma \ref{lemA.1}, to prove Lemma \ref{lemA.2}, we only need to consider the following two bias terms
$$\sup_{\textbf{u}\in\mathfrak{X},\theta\in\Theta}
\Big|f_{\theta^T\textbf{X}}(\theta^T\textbf{u})-\E\tilde{f}_{\theta^T\textbf{X},n}(\theta^T\textbf{u})\Big|$$
and
$$\sup_{\textbf{u}\in\mathfrak{X},\theta\in\Theta}
\Big|\phi_{\theta^T\textbf{X}}(\theta^T\textbf{u})-\E\tilde{\phi}_{\theta^T\textbf{X},n}(\theta^T\textbf{u})\Big|.$$
From the classic change of variable,
a Taylor expansion and Conditions (C2) and (C6),
we  get that both the bias terms are
of order $O(h^2)$.
Hence, we complete the proof of Lemma \ref{lemA.2}.
\qed

\begin{lem}\label{lemA.3}
Under Conditions (C1)-(C3), (C5), (C6) and (C8), we have
\begin{eqnarray}\label{EqA.7}
&&\sup_{\textbf{u}\in\mathfrak{X},\theta\in\Theta}
\Big|f_{\theta^T\textbf{X}}(\theta^T\textbf{u})-\hat{f}_{\theta^T\textbf{X},n}(\theta^T\textbf{u})\Big|=
O\Big(\sqrt{\frac{\log n}{nh}}+h^2\Big)~a.s.
\end{eqnarray}
and
\begin{eqnarray}\label{EqA.8}
&&\sup_{\textbf{u}\in\mathfrak{X},\theta\in\Theta}
\Big|\phi_{\theta^T\textbf{X}}(\theta^T\textbf{u})-\hat{\phi}_{\theta^T\textbf{X},n}(\theta^T\textbf{u})\Big|=
O\Big(\sqrt{\frac{\log n}{nh}}+h^2\Big)~a.s.
\end{eqnarray}
\end{lem}
\noindent{\it Proof:}
We first consider \eqref{EqA.7}.
Similar to the proof of  Lemma 2 in  Lemdani et al. \cite{LOP2009}, we get
from Theorem 3.2 in He and Yang \cite{HY1998b},
Theorem 4.1 in He and Yang \cite{HY1998a} and the strong law of large numbers that
$$\sup_{\textbf{u}\in\mathfrak{X},\theta\in\Theta}\Big|\tilde{f}_{\theta^T\textbf{X},n}
(\theta^T\textbf{u})-\hat{f}_{\theta^T\textbf{X},n}(\theta^T\textbf{u})\Big|=O(n^{-1/2})~a.s.,$$
which, together with Lemma \ref{lemA.2} and Condition (C8), implies that \eqref{EqA.7} holds.

Following the same lines as  the proof of \eqref{EqA.7}, we  get \eqref{EqA.8} by
using Lemma \ref{lemA.2} again.
\qed

Noting that
$g(s;\theta)=\frac{\phi_{\theta^T\textbf{X}}(s)}{f_{\theta^T\textbf{X}}(s)}$
and
$\hat{g}_{n}(s;\theta)=\frac{\hat{\phi}_{\theta^T\textbf{X},n}(s)}{\hat{f}_{\theta^T\textbf{X},n}(s)}$,
we get the following lemmas from Lemma \ref{lemA.3}.
\begin{lem}\label{lemA.4}
Under the assumptions of Lemma \ref{lemA.3}, we have
\begin{eqnarray*}
&&\sup_{\textbf{u}\in\mathfrak{X},\theta\in\Theta}
\Big|g(\theta^T\textbf{u};\theta)-\hat{g}_{n}(\theta^T\textbf{u};\theta)
\Big|=O\Big(\sqrt{\frac{\log n}{nh}}+h^2\Big)~a.s.
\end{eqnarray*}
\end{lem}

Similar to the proof of Lemma \ref{lemA.4}, we have
\begin{lem}\label{lemA.5}
Under the assumptions of Lemma \ref{lemA.3}, we have
\begin{eqnarray*}
&&\sup_{\textbf{u}\in\mathfrak{X},\theta\in\Theta}
\Big|\nabla_{\theta}g(\theta^T\textbf{u};\theta)
-\nabla_{\theta}\hat{g}_{n}(\theta^T\textbf{u};\theta)\Big|=
O\Big(\sqrt{\frac{\log n}{nh^3}}+h^2\Big)~a.s.
\end{eqnarray*}
\end{lem}

\subsection{Proofs of Theorems \ref{Th3.1} and \ref{Th3.2}}

In this subsection, we give the detailed proofs of Theorems \ref{Th3.1} and \ref{Th3.2}.
Define
\begin{eqnarray*}
&&M(\theta,g)=\int [v-g(\theta^T\textbf{u})]^2J(\textbf{u})H(d\textbf{u},dv).
\end{eqnarray*}

\noindent{\it Proof of Theorem \ref{Th3.1}:}
From Theorem 5.7 in van der Vaart \cite {V1998}, to prove Theorem \ref{Th3.1},
we only need to show that
\begin{equation}\label{EqA.9}
\sup_{\theta\in\Theta}|M_n(\theta,\hat{g}_n)-M(\theta,g)|=o_\P(1),
\end{equation}
where $M_n(\theta,g)$ is defined in \eqref{Eq2.3}. 
We first consider the difference
\begin{eqnarray*}
&&|M_n(\theta,\hat{g}_n)-M_n(\theta,g)|\\
&&\quad\leq
\sup_{\textbf{u}\in\mathfrak{X},\theta\in\Theta}\Big|\hat{g}_n(\theta^T\textbf{u})-g(\theta^T\textbf{u})\Big|
\int\Big[\sup_{\textbf{u}\in\mathfrak{X},\theta\in\Theta}|g(\theta^T\textbf{u})+\hat{g}_n(\theta^T\textbf{u})|+2|v|\Big]H_n(d\textbf{u},dv).
\end{eqnarray*}
From Conditions (C2) and (C4), the integral on the righthand side is finite.
By Lemma \ref{lemA.4}, we  deduce that
$$\sup_{\theta\in\Theta}|M_n(\theta,\hat{g}_n)-M_n(\theta,g)|=o_\P(1).$$
Moreover, similar to the proof of Theorem 1.1 in Stute \cite{S1999},
we get $$\sup_{\theta\in\Theta}|M_n(\theta,g)-M(\theta,g)|=o_\P(1).$$
Hence \eqref{EqA.9} holds and we end the proof of Theorem \ref{Th3.1}.
\qed

To get the asymptotic normality of $\hat{\theta}_n$,
we first consider the case $g$ is known.
From Sherman \cite[Theorems 1, 2]{S1994}, to prove our result,
we only need to study the representation of $M_n(\theta,g)$. In fact, we have the following lemma.

\begin{lem}\label{lemA.6} \quad
\begin{itemize}
\item[(i)] Under  Conditions (C1)-(C4), we have, on $o_\P(1)$
neighborhoods of $\theta_0$,
\begin{equation}\label{EqA.10}
M_n(\theta,g)=M(\theta,g)+O_\P\Big(\frac{\|\theta-\theta_0\|}{\sqrt{n}}\Big)+o_\P(\|\theta-\theta_0\|^2)
+K_n(\theta_0),
\end{equation}
where
 \begin{eqnarray*}
K_n(\theta_0)=\int[v-g(\theta_0^T\textbf{u})]^2
J(\textbf{u})\Big(H_n(d\textbf{u},dv)-H(d\textbf{u},dv)\Big).
\end{eqnarray*}

\item[(ii)] Under Conditions (C1)-(C4), we have, on $O_{\P}(n^{-1/2})$
neighborhoods of $\theta_0$,
\begin{equation}\label{EqA.11}
M_n(\theta,g)=\frac{1}{2}(\theta-\theta_0)^T\Lambda(\theta-\theta_0)
+n^{-1/2}(\theta-\theta_0)^TW_n+o_{\P}(n^{-1})+K_n(\theta_0),
\end{equation}
where
$W_n=n^{-1/2}\sum_{i=1}^n\zeta_i(\psi)$ is a random vector,
$\psi(\cdot)$ and $\zeta_i(\cdot)$ are defined by \eqref{Eq3.3} and
 \eqref{Eq3.4}, respectively.
 \end{itemize}
\end{lem}

\noindent{\it Proof:}
We only need to show \eqref{EqA.10}. \eqref{EqA.11} can be done in the same way.
Note that
\begin{eqnarray*}
M_n(\theta,g)-M(\theta,g)
&&=2\int[v-g(\theta_0^T\textbf{u})]\big[g(\theta_0^T\textbf{u})-g(\theta^T\textbf{u})\big]
J(\textbf{u})\Big(H_n(d\textbf{u},dv)-H(d\textbf{u},dv)\Big)\nonumber\\
&&\quad\quad+\int\big[g(\theta_0^T\textbf{u})-g(\theta^T\textbf{u})\big]^2
J(\textbf{u})\Big(H_n(d\textbf{u},dv)-H(d\textbf{u},dv)\Big)\nonumber\\
&&\quad\quad+\int[v-g(\theta_0^T\textbf{u})]^2J(\textbf{u})\Big(H_n(d\textbf{u},dv)-H(d\textbf{u},dv)\Big)\\
&&=:A_{1n}+A_{2n}+K_n(\theta_0).
\end{eqnarray*}
From a Taylor's expansion, $A_{1n}$ can be rewritten as
\begin{eqnarray}\label{EqA.12}
&&2(\theta_0-\theta)^T\int\psi(\textbf{u},v)\Big(H_n(d\textbf{u},dv)-H(d\textbf{u},dv)\Big)\nonumber\\
&&\quad\quad-(\theta_0-\theta)^T\int \beta(\textbf{u},v)\Big(H_n(d\textbf{u},dv)-H(d\textbf{u},dv)\Big)(\theta_0-\theta),
\end{eqnarray}
where $\psi(\textbf{u},v)$ is defined by \eqref{Eq3.3},
and
$$\beta(\textbf{u},v)=[v-g(\theta_0^T\textbf{u})]\nabla_{\theta,\theta}^2g(\theta_1^T\textbf{u})
J(\textbf{u})$$
with $\theta_1$ being a vector between $\theta$ and $\theta_0$.
It follows from Conditions (C1) to (C4) that \eqref{EqA.1} holds for $\psi(\textbf{u},v)$.
Hence, by applying Corollary \ref{corA.1}, the first term in \eqref{EqA.12} is
$$2(\theta_0-\theta)^T\Big[n^{-1}\sum_{i=1}^n\zeta_i(\psi)+
o_{\P}(n^{-\frac{1}{2}})\Big].$$
It follows from the multivariate central limit theorem
that the first term in \eqref{EqA.12} is of order $O_\P\Big(\frac{\|\theta-\theta_0\|}{\sqrt{n}}\Big)$.
By the strong consistency of the Lynden-Bell integral (He and Yang \cite[Theorem 3.2]{HY2003}) and
the boundedness of $\nabla_{\theta}g$ and $\nabla^2_{\theta,\theta}g$ (see Condition (C4)),
 the second term in \eqref{EqA.12} is $o_\P(\|\theta-\theta_0\|^2)$.
Moreover, by a  Taylor's expansion,
\begin{eqnarray*}
A_{2n}=(\theta_0-\theta)^T\int\Big[\nabla_\theta g(\theta_2^T\textbf{u})
\nabla_\theta g(\theta_2^T\textbf{u})^T\Big]
J(\textbf{u})\Big(H_n(d\textbf{u},dv)-H(d\textbf{u},dv)\Big)(\theta_0-\theta),
\end{eqnarray*}
where $\theta_2$ is a vector between $\theta$ and $\theta_0$.
Using the strong consistency of the Lynden-Bell integral and Condition (C4) again,
we obtain that $A_{2n}$ is also of order $o_\P(\|\theta-\theta_0\|^2)$.
The proof of \eqref{EqA.10} is completed.
\qed

With Lemmas \ref{lemA.4} to \ref{lemA.6} in hand,
we are able to present the proof of
 Theorem \ref{Th3.2} in the following.

\noindent{\it Proof of Theorem \ref{Th3.2}:}
The proof of Theorem \ref{Th3.2} is similar to that of the Main Lemma
in Bouaziz and Lopez \cite{BL2010}.
From Theorems 1 and 2 of Sherman \cite{S1994} and Lemma \ref{lemA.6}, to prove Theorem \ref{Th3.2}, we only need to show that
\begin{equation}\label{EqA.13}
M_n(\theta,\hat{g}_n)=M_n(\theta,g)+\tilde{K}_{n}(\theta_0)+
o_\P\Big(\frac{\|\theta-\theta_0\|}{\sqrt{n}}\Big)+o_\P(\|\theta-\theta_0\|^2),
\end{equation}
where $\tilde{K}_{n}(\theta_0)$ is a term that depends only on $\theta_0$.
Following similar ideas as those in the proof of Theorem 3.5 in \cite{LPV2013},
especially the theory of empirical process,
we obtain that \eqref{EqA.13} holds. We end the proof of Theorem \ref{Th3.2}.
\qed

\end{appendix}

\noindent{\bf Acknowledgments:}\ Dr. Kong was supported by the National Natural
Science Foundation of China (Nos. 11601260 and 71671104),
 the Project of Humanities and Social Science of Ministry of Education of China (No.16YJA910003).
 Dr. Zhang was supported by Project of Shandong Provincial Higher Educational
Science and Technology Program (Nos. J16LI56 and J17KA163),
the Fostering Project of Dominant Discipline and Talent Team of Shandong University of Finance and Economics.
 Dr. Dai was supported by the National Natural Science Foundation of China (No.11361007) and
 the Fostering Project of Dominant Discipline and Talent Team of Shandong Province Higher Education Institutions.


\begin{thebibliography}{99}

\bibitem{BFZ2009}
Bai, Y., Fung, W., Zhu, Z., 2009.
Penalized quadratic inference functions
for single-index models with longtitudinal data.
\emph{J. Multivariate Anal.} {\bf 100}, 152-161.

\bibitem{BL2010}
Bouaziz, O., Lopez, O., 2010.
Conditional density estimation in a censored single-index regression model.
\emph{Bernoulli} {\bf 16(2)},  514-542.


\bibitem{CRA2006}
Chaieb, L.L., Rivest, L.P., Abdous, B., 2006.
 Estimating survival under a dependent
truncation. \emph{Biometrika} {\bf 93}, 655-669.

\bibitem{DHP2006}
Delecroix, M., Hristache, M., Patilea, V., 2006.
 On semiparametric M-estimation in single-index regression.
 \emph{J. Statist. Plann. Inference} {\bf 136}, 730-769.


\bibitem{DL1991}
Duan, N., Li, K.C., 1991.
Slicing regression: a link-free regression method.
\emph{Ann. Statist.} {\bf 19}, 505-530.

 \bibitem{EM2000}
Einmahl, U., Mason, D.M., 2000.
An empirical process approach to the uniform consistency of kernel-type function estimators.
\emph{J. Theoret. Probab.} {\bf 13}, 1-37.

\bibitem{GG1999}
Gin\'{e}, E., Guillou, A., 1999.
Laws of the iterated logarithm for censored data. \emph{Ann. Probab.} {\bf 27},
 2042-2067.

\bibitem{HT1993}
H\"{a}rdle, W., Tsybakov, A.B., 1993.
How sensitive are average derivatives?
\emph{J. Econometrics} {\bf 58}, 31-48.

\bibitem{HHI1993}
H\"{a}rdle, W., Hall, P., Ichimura, H., 1993.
Optimal smoothing in single-index models. \emph{Ann. Statist.}  {\bf  21}, 157-178.

\bibitem{HY1998a}
He, S., Yang, G.L., 1998.
The strong law under random truncation.
\emph{Ann. Statist.} {\bf 26},
992-1010.

\bibitem{HY1998b}
He, S., Yang, G.L., 1998.
Estimation of the truncation probability in the random truncation model.
\emph{Ann. Statist.} {\bf 26}, 1011-1027.

\bibitem{HY2003}
He, S., Yang, G.L., 2003.
Estimation of regression parameters with left truncated data.
\emph{J. Statist. Plann. Inference}
{\bf 117}, 99-122.

\bibitem{I1993}
Ichimura, H., 1993.
Semiparametric least squares (SLS) and weighted SLS estimation of single-index models.
\emph{J. Econometrics} {\bf 58}, 71-120.

\bibitem{KX2007}
Kong, E.F., Xia, Y.C., 2007.
Variable selection for the single index model.
\emph{Biometrika} {\bf 941}, 217-229.

\bibitem{LOP2009}
Lemdani, M., Ould-Sa\"{\i}d, E., Poulin, N., 2009.
Asymptotic properties of a conditional quantile estimator with randomly truncated data.
\emph{J. Multi. Analy.} {\bf 100}, 546-559.

\bibitem{L2009}
Lopez, O., 2009.
Single-index regression models with right-censored responses.
\emph{J. Statist. Plann. Inference} {\bf 139}, 1082-1097.

\bibitem{LPV2013}
Lopez, O., Patilea, V., Van Keilegom, I., 2013.
Single index regression models in the presense of
censoring depending on the covariates.
\emph{Bernoulli} {\bf 19}, 721-747.

\bibitem{LB2005}
Lu, X., Burke, M.D., 2005.
Censored multiple regression by the method of average derivatives.
\emph{J. Multivariate Anal.} {\bf 95}, 182-205.

\bibitem{L1971}
Lynden-Bell, D., 1971.
A method of allowing for known observational selection in small samples
applied to 3CR quasars.
 \emph{Monthly Notices Roy. Astronom. Soc.} {\bf 155}, 95-118.

 \bibitem{MUM2016}
Moreira, C., de U\~{n}a-\'{A}lvarez, J., Meira-Machado, L., 2016.
Nonparametric regression with doubly truncated data.
 \emph{Comput. Statist. Data Anal.} {\bf 93}, 294-307.

\bibitem{OL2006}
Ould-Sa\"{\i}d, E., Lemdani, M., 2006.
Asymptotic properties of a nonparametric
regression function estimator with randomly truncated data.
\emph{Ann. Instit. Statist. Math.} {\bf 58}, 357-378.


\bibitem{S2006}
S\'{a}nchez Sellero, C., Gonz\'{a}lez Manteiga, W., Van Keilegom, I., 2005.
Uniform representation
of product-limit integrals with applications.
\emph{Scand. J. Statist. } {\bf 32}, 563-581.

\bibitem{S1980}
Serfling, R.J., 1980.
\emph{Approximation Theorems of Mathematical Statistics.}
Wiley: New York.

\bibitem{S1994}
Sherman, R.P., 1994.
Maximal inequalities for degenerate U-processes with applications to optimization estimators.
\emph{Ann. Statist.} {\bf 22}, 439-459.

\bibitem{S1986}
Stoker, T.M., 1986.
Consistent estimation of scaled coefficients.
\emph{Econometrica } {\bf 54}, 1461-1481.

\bibitem{S1993}
Stute, W., 1993.
Almost sure representations of the product-limit estimator for truncated data.
 \emph{Ann. Statist.} {\bf 21}, 146-156.

\bibitem{S1996}
Stute, W., 1996.
Distributional convergence under random censorship when covariables are present.
 \emph{Scan. J. Statist.} {\bf 23}, 461-471.

\bibitem{S1999}
Stute, W., 1999.
Nonlinear censored regression.
 \emph{Statist. Sinica} {\bf 9},
1089-1102.


\bibitem{SW2008}
Stute, W., Wang, J.L., 2008.
The central limit theorem under
random truncation.
\emph{Bernoulli} {\bf 14},
604-622.

\bibitem{V1998}
van der Vaart, A.W., 1996.
\emph{Asymptotics Statistics.
Cambridge Series in Statistical and Probabilistic Mathematics} {\bf 3.}
Cambridge Univ. Press: Cambridge.

\bibitem{VW1996}
van der Vaart, A.W., Wellner, J.A., 1996.
\emph{Weak Convergence and Empirical Processes.}
Springer-Verlag: New York.

\bibitem{W2009}
Wang, H.B., 2009.
Bayesian estimation and variable selection for single index models.
\emph{Comput. Statist. Data  Anal.} {\bf 53}, 2617-2627.

\bibitem{WSHW2010}
Wang, Y., Shen, J., He, S., Wang, Q., 2010.
Estimation of single index model with missing response at random.
\emph{J. Statist. Plann. Inference} {\bf 140}, 1671-1690.

\bibitem{WY2009}
Wang, L., Yang, L., 2009.
Spline estimation of  single-index models.
\emph{Statist. Sinica},
{\bf 19}, 765-783.


\bibitem{W1985}
Woodroofe, M., 1985.
Estimating a distribution function with truncated data.
\emph{Ann. Statist.} {\bf 13},
163-177.

\bibitem{XLTZ2009}
Xia, Y.C., Li, W.K., Tong, H., Zhang, D.X., 2004.
A goodness-of-fit test for single-index models.
\emph{Statist. Sinica} {\bf 14}, 1-39.

\bibitem{XZ2006}
Xue, L.G., Zhu, L.X., 2006.
Empirical likelihood for single-index models.
\emph{J. Multi. Analy.} {\bf 97}, 1295-1312.

\bibitem{YC2005}
Yin, X.R., Cook, R.D., 2005.
Direction estimation in single-index regressions.
\emph{Biometrika} {\bf 92}, 371-384.

\bibitem{Z2011}
Zhou, W., 2011.
A weighted quantile regression for randomly truncated data.
\emph{Comput. Statist. Data Anal.}
{\bf 55}, 554-566.

\end{thebibliography}
\end{document}